\documentclass[reqno,11pt]{amsart}

\usepackage[all]{xy}
\usepackage{hhline}
\usepackage{amssymb}
\usepackage{amsmath,amsthm,amscd,amsfonts}

\usepackage[dvips, dvipsnames, usenames]{color}

\usepackage[active]{srcltx}

\def\k{\Bbbk}

\newcommand{\GL}{\mathbf{GL}}

\newcommand{\ot}{\otimes}

\newcommand{\com}{\Delta}

\newcommand\toba{{\mathfrak B }}

\newcommand{\gr}{\operatorname{gr}}
\newcommand{\trid}{\triangleright}

\newcommand{\Lc}{{\mathcal L}}

\newcommand{\DD}{{\mathbb D}}

\newcommand{\K}{{\mathcal K}}
\newcommand{\Z}{{\mathbb Z}}

\newcommand{\D}{{\mathcal D}}

\newcommand{\II}{\mathcal{I}}

\newcommand{\Oc}{{\mathcal O}}
\newcommand{\oc}{{\mathcal O}}
\newcommand{\ydh}{{}^H_H\mathcal{YD}}

\newcommand{\supp}{\operatorname{supp}}

\newcommand{\End}{\operatorname{End}}
\newcommand{\Aut}{\operatorname{Aut}}

\newcommand\card{\operatorname{card}}

\newcommand\sgn{\operatorname{sgn}}

\newcommand\co{\operatorname{co}}
\newcommand{\ydho}{{}^{H_0}_{H_0}\mathcal{YD}}

\theoremstyle{plain}

\newtheorem{lema}{Lemma}[section]
\newtheorem{theorem}[lema]{Theorem}
\newtheorem{cor}[lema]{Corollary}

\newtheorem{prop}[lema]{Proposition}

\theoremstyle{definition}
\newtheorem{definition}[lema]{Definition}
\newtheorem{exa}[lema]{Example}

\theoremstyle{remark}
\newtheorem{obs}[lema]{Remark}
\newtheorem{rmk}[lema]{Remarks}

\theoremstyle{plain}
\newcounter{maint}

\newtheorem{mainthm}[maint]{Theorem}

\theoremstyle{plain}

\newcommand\id{\operatorname{id}}

\newcommand\dm{\mathbb D_m}

\newcommand\s{\mathbb S}

\newcommand{\yddm}{{}^{\k \dm}_{\k \dm}\mathcal{YD}}

\def\pf{\begin{proof}}
\def\epf{\end{proof}}

\theoremstyle{remark}

\begin{document}

\renewcommand{\baselinestretch}{1.2}

\thispagestyle{empty}

\title[On pointed Hopf algebras over dihedral groups]
{On pointed Hopf algebras over dihedral groups}

\author[Fernando Fantino, Gast\'on A. Garc\'ia]{Fernando Fantino
and Gast\'on Andr\'es Garc\'ia}

\thanks{This work was partially supported by
Universidad Nacional de C\'ordoba,
ANPCyT-Foncyt, CONICET,
Ministerio de Ciencia y Tecnolog\'\i a (C\'ordoba) and Secyt (UNC)}

\address{\newline\noindent Facultad de Ciencias Exactas, F\'isicas y Naturales \&
\newline\noindent Facultad de Matem\'atica, Astronom\'\i a y F\'\i sica.
\newline\noindent
Universidad Nacional de C\'ordoba.
\newline \noindent
CIEM -- CONICET.
\newline \noindent
Medina Allende s/n, Ciudad Universitaria
\newline \noindent 5000 C\'ordoba,
Argentina}

\email{(fantino|ggarcia)@famaf.unc.edu.ar}


\begin{abstract}
Let $\k$ be an algebraically closed field of characteristic $0$ and
let $\dm$ be the dihedral group of order $2m$ with
$m= 4t, t\geq 3$.
We classify all finite-dimensional
Nichols algebras over $\dm$ and all
finite-dimensional pointed Hopf
algebras whose group of group-likes is $\dm$,
by means of the lifting method. As a byproduct
we obtain new examples of finite-dimensional
pointed Hopf algebras.
\end{abstract}

\maketitle


\section*{Introduction}
This paper is concerned with the classification
of finite-dimensional Hopf algebras over an
algebraically closed field $\k$ of characteristic $0$.
In particular, we study pointed Hopf algebras over
dihedral groups $\dm$, $m=4t\geq 12$, using
the lifting method, which leads to the study of
finite-dimensional Nichols algebras in the
category $\yddm$ of left Yetter-Drinfeld modules over
$\dm$. For more examples over $\mathbb{D}_{p}$, 
$p$ an odd prime or $4$, see \cite[Section 3.3]{AG1}.

\smallbreak A significant progress has been achieved in \cite{AS3}
in the case of pointed Hopf algebras with abelian group of
group-likes. When the group of group-likes is not abelian, the
problem is far from being completed. Some hope is present in the
lack of examples: in this situation, Nichols algebras tend to be
infinite dimensional, see for example \cite{AZ,AF2,AFZ,AFGV, FGV,
FGV2}. Nevertheless, examples on which the Nichols algebras are
finite dimensional do exist. Over $\s_{3}$ and $\s_4$ these
algebras were determined in \cite{AHS}. All of them arise from
racks associated to a cocycle, and in \emph{loc. cit} and
\cite{GG} the classification of pointed Hopf algebras over
$\s_{3}$ and $\s_4$ is completed, respectively.

\par Let $G$ be a finite group and let $A_0$
be the group algebra of
$G$. The main steps of the lifting method for
the classification of all finite-dimensional pointed Hopf algebras
with group $G$ are:

\begin{enumerate}
\item[$(a)$] determine all Yetter-Drinfeld modules $V$ such that the
Nichols algebra $\toba(V)$ is finite dimensional,

\item[$(b)$] for such $V$, compute all Hopf algebras $A$ such
that $\gr A\simeq \toba(V)\# A_{0}$, the Radford-Majid product. 
We call $A$ a \emph{lifting}
of $\toba(V)$ over $G$.

\item[$(c)$] Prove that any finite-dimensional pointed
Hopf algebras with group $G$
is generated by group-likes and skew-primitives.
\end{enumerate}

Assume $G=\dm$, $m=4t, t\geq 3$.
In Section \ref{sec:nichols-dihedral} we complete step $(a)$,
that is, we determine all $V \in \yddm$ such that the Nichols
algebra $\toba(V)$ is finite-dimensional, and we
describe explicitly these Nichols algebras.
Then we prove step $(b)$ and $(c)$ in
Section \ref{sec:pointed-ha-dm}, which are
given by Theorem \ref{teo:liftings-quadratic}
and Theorem \ref{teo:item-c}, respectively.

\par Summarizing, the main
theorems of the present paper are the
following; for definitions see Definitions \ref{def:I}, \ref{def:L} and
\ref{def:I-L}.
\begin{mainthm}\label{teo:all-nichols-dm}
Let $\toba(M)$ be a finite-dimensional Nichols algebra
in $ \yddm $. Then $\toba(M)\simeq \bigwedge M$, with $M$
isomorphic
either to $M_{I}$, or to
$M_{L}$, or to $M_{I,L}$, with $I\in \II$, $L\in \Lc$ and  
$(I,L) \in \K$, respectively.
\end{mainthm}

The proof of the preceeding theorem uses the 
classification of finite-dimensional Nichols algebras 
of diagonal type due to
I. Heckenberger \cite{H2}. 

Although all Nichols algebras in $\yddm$ turn out
to be exterior algebras, we write $\toba(M)$ to enphasize
the Yetter-Drinfeld module structure.
The following theorem gives all liftings of
these families of Nichols algebras.

\begin{mainthm}\label{thm:B}\label{teo:liftings-quadratic}
\label{prop:quad-are-liftings}
Let $H$ be a finite-dimensional pointed
Hopf algebra over $\DD_{m}$.
Then $H$ is isomorphic to one of the following algebras
\begin{enumerate}
\item[$ (a) $] $\toba(M_I)\#\k \DD_{m}$, with $I= \{(i,k)\} \in \II$, $k\neq n$.
\item[$ (b) $] $\toba(M_{L})\#\k \DD_{m}$, with $L\in \Lc$.
\item[$ (c) $] $A_{I}(\lambda,\gamma)$, with $I\in \II$, $|I| > 1$
or $I = \{(i,n)\}$ and $\gamma \equiv 0$.
\item[$ (d) $] $B_{I,L }(\lambda,\gamma,\theta,\mu)$, with 
$(I,L)\in \K$, $|I|> 0$ and $|L|>0$.
\end{enumerate}
Conversely, any Hopf algebra appearing in the list
is a lifting
of a finite-dimensional Nichols algebra in $ \yddm $.
\end{mainthm}

After the study of finite-dimensional
pointed Hopf algebras over $\s_{3}, \s_{4}$,  
this theorem is the first result that 
gives an infinite family of non-abelian
groups where the classification of finite-dimensional pointed Hopf 
algebras with non-trivial examples is completed and, unlike the symmetric groups case,
it provides for each dihedral group infinitely many non-trivial
finite-dimensional pointed Hopf algebras.  

The paper is organized as follows. In Section \ref{sec:prelim}
we establish conventions and recall some basic
facts about pointed Hopf algebras $H$ such as
the coradical filtration, the grading associated to it and
the category $\ydho $ of Yetter-Drinfeld modules over the corradical.
If $G=G(H)$, the irreducible modules of  $\ydho$
are parametrized by pairs $(\oc,\rho)$, where $\oc$
is a conjugacy class of $G$ and $\rho$ is a simple
representation of the centralizer of an element $\sigma \in \oc$.
At the end of this first section we recall the \textit{type D}-criterium
\cite[Thm. 3.6]{AFGV}, which
helps to determine when the Nichols algebra $\toba(\oc,\rho)$
associated to $(\oc, \rho)$ is infinite-dimensional,
depending only on the rack
structure of the conjugacy class $\oc$.
In Section \ref{sec:nichols-dihedral} we work with Nichols
algebras over the dihedral groups $\dm$, with $m=4t\geq 12$ and
give the proof of Theorem \ref{teo:all-nichols-dm}.
We begin by determining which irreducible modules give rise
to finite-dimensional Nichols algebras and then we
extend our study to arbitrary modules. It turns out
that all finite-dimensional Nichols algebras in $\yddm$
are exterior algebras of some irreducible modules or
specific families of them.
The last section of the paper is devoted to the classification
of pointed Hopf algebras over $\dm$, that is, to the proof of Theorem
\ref{thm:B}. It consists mainly in the construction of
the liftings of the finite-dimensional Nichols algebras
given in Section \ref{sec:nichols-dihedral}. To do this, we
show first in Theorem \ref{teo:item-c}
that all pointed Hopf algebras over $\dm$, $m=4t\geq 12$
are generated by group-likes and skew-primitive elements.
Then we prove that if $H$ is a pointed Hopf algebra
over $\dm$, then some quadratic relations must hold and
using these relations we define in Definitions \ref{liftings1} and
\ref{def:liftings2} two families of quadratic algebras. Finally,
using representation theory we prove that these algebras together
with the bosonizations are
all the possible liftings.
We conclude the paper with the study of the isomorphism classes.

\section{Preliminaries}\label{sec:prelim}

\subsection{Conventions}We
work over an algebraically closed field $\k$ of characteristic
zero. Let $H$ be a Hopf algebra over $\k$ with bi\-ject\-ive
antipode. We use Sweedler's notation $\Delta(h)=h_1\ot
h_2$ for the comultiplication in $H$, but dropping the summation
symbol, see \cite{S}.

The \emph{coradical}
$H_0$ of $H$ is the sum of all simple
sub-coalgebras of $H$. In particular, if $G(H)$ denotes the group
of \emph{group-like elements} of $H$, we have $\k G(H)\subseteq
H_0$. We say that a Hopf algebra is \emph{pointed} if $H_0=\k
G(H)$. Denote by $\{H_i\}_{i\geq 0}$ the \emph{coradical
filtration} of $H$; if $H_0$ is a Hopf
subalgebra of $H$, then $\gr H = \oplus_{n\ge 0}\gr H(n)$ is the associated
graded Hopf algebra, with
 $\gr H(n) = H_n/H_{n-1}$ (set
$H_{-1} =0$).  Let $\pi:\gr H\to H_0$ be the homogeneous
projection, then $R= (\gr H)^{\co \pi}$ is the \emph{diagram} of
$H$; which is a braided Hopf algebra in the category $\ydho$
of left Yetter-Drinfeld modules over $H_{0}$, and it is a
graded sub-object of $\gr H$. The linear space
$R(1)$, with
the braiding from $\ydho$, is called the \emph{infinitesimal
braiding} of $H$ and coincides with the subspace of primitive
elements $P(R)= \{r \in R:\ \Delta_{R}(r) = r\ot 1 + 1\ot r\}$.
It turns out that the Hopf algebra $\gr H$ is the
Radford-Majid biproduct $\gr H \simeq R\# \k G(H)$ and the
subalgebra of $R$ generated by $V$ is isomorphic to the Nichols
algebra $\toba(V)$.

\subsection{Yetter-Drinfeld modules over $\k G$}\label{subsect:ydCG}
Let $G$ be a finite group. A left Yetter-Drinfeld module over $\k G$ is a left
$G$-module and left $\k G$-comodule $M$ such that
$$
\delta(g.m) = ghg^{-1} \otimes g.m, \qquad \forall\ m\in M_h, g, h\in G,
$$
where $M_h = \{m\in M: \delta(m) = h\otimes m\}$; clearly, $M =
\oplus_{h\in G} M_h$. The \emph{support} of $M$ is $\supp
M = \{g\in G: M_g \neq 0\}$. Yetter-Drinfeld modules over $G$ are
completely reducible. Also, irreducible Yetter-Drinfeld modules
over $G$ are parameterized by pairs $(\oc, \rho)$, where $\oc$ is a
conjugacy class and $(\rho,V)$ is an irreducible representation of the
centralizer $C_G(\sigma)$ of a fixed point $\sigma\in \oc$. We
denote the corresponding Yetter-Drinfeld module by
$M(\oc, \rho)$ and by $\toba(\oc,\rho)$ the associated
Nichols algebra.

Here is a precise description of the Yetter-Drinfeld module $M(\oc,
\rho)$. Let $\sigma_1 = \sigma$, \dots, $\sigma_{n}$
be a numeration of $\oc$ and
let $g_i\in G$ such that $g_i \sigma g_{i}^{-1} = \sigma_i$ for all $1\le i \le
n$. Then  $M(\oc, \rho) = \oplus_{1\le i \le n}g_i\otimes V$. Let
$g_iv := g_i\otimes v \in M(\oc,\rho)$, $1\le i \le n$, $v\in V$.
If $v\in V$ and $1\le i \le n$, then the action of $g\in G$ and
the coaction are given by
$$g\cdot (g_iv) = g_j(\gamma\cdot v), \qquad
\delta(g_iv) = \sigma_i\otimes g_iv,
$$
where $gg_i = g_j\gamma$, for some $1\le j \le n$ and $\gamma\in
C_G(\sigma)$. The explicit formula for the braiding is then given by
\begin{equation} \label{yd-braiding}
c(g_iv\otimes g_jw) = \sigma_i\cdot(g_jw)\otimes g_iv = g_h(\gamma\cdot
v)\otimes g_iv\end{equation} for any $1\le i,j\le n$, $v,w\in V$,
where $\sigma_ig_j = g_h\gamma$ for unique $h$, $1\le h \le n$ and
$\gamma \in C_G(\sigma)$.
Since
$\sigma\in Z(C_G(\sigma))$, Schur's Lemma says that
$\sigma$  acts by a scalar $q_{\sigma\sigma}$
on $V$.

\medbreak

The following are useful tools that, under certain conditions,
allow us to determine if the dimension of a Nichols algebra is
infinite. These results are about abelian and non-abelian subracks
of a conjugacy class $\oc$ of $G$, respectively.

\begin{lema}\label{le:odd}\cite[Lemma
2.2]{AZ} Let $G$ be a finite group, $\oc_{\sigma}$ a conjugacy
class in $G$. If $\oc_{\sigma}$ is real (\textit{i.e.}
$\sigma^{-1}\in \Oc$)
 and $\dim\toba(\oc_{\sigma}, \rho)< \infty$,
then $q_{\sigma\sigma} = -1$ and $\sigma$ has even order.\qed
\end{lema}

We say that $\oc$
 is \emph{of type D} if there exist $r$, $s\in \oc$ such that
 $(rs)^2\neq (sr)^2$ and
$r$ and $s$ are not conjugate in some subgroup
$H$ of $G$ containing $r$ and $s$.

\begin{lema}\cite[Thm. 3.6]{AFGV}.\label{lem:tipo-D-grupo}
If $\oc$ is of type $D$, then $\toba(\oc, \rho)$ is
infinite-dimensional for all $\rho$.\qed
 \end{lema}

Let $ A $ be a finite abelian group and $g\in \Aut (A)$.
We denote by $(A,g)$ the rack with underlying set $A$ and 
rack multiplication $x\trid y := g(y) + (\id-g)(x)$, $x, y\in A$;
this is a subrack of the group $A \rtimes \langle g\rangle$.
Any rack isomorphic to some $(A,g)$ is called \emph{affine}.

For instance, consider the cyclic group $A = C_n$ and 
the automorphism $g$ given by the inversion;
the rack $(A,g)$ is denoted $\D_{n}$ and called a dihedral rack.
Thus, a family
$(\mu_{i})_{i\in \Z/n}$ of distinct elements of a rack $X$
is isomorphic to $\D_{n}$ if $\mu_{i} \trid \mu_{j} = \mu_{2i-j}$ for all $i,j$.

\begin{lema}\cite[Lemma 2.1]{thr}.\label{lema:d2m}
If $m>2$, then the dihedral rack $\D_{2m}$ is of type D.\qed
\end{lema}

\section{Nichols algebras over $\dm, m=4t\geq 12$}
\label{sec:nichols-dihedral}

Let $m$ be a positive integer, $m\geq 3$. The dihedral group of
order $2m$ can be presented by generators and relations as follows
\begin{align}\nonumber
\mathbb D_m&:=\langle x,y \quad | \quad x^2=1=y^m  \,\, , \,\,
xy=y^{-1}x \rangle.
\end{align}

\textit{From now on we assume that $m=4t$, with $t\geq 3$, and
set $n=\frac{m}{2}= 2t$.}

\bigbreak In this section we determine all finite-dimensional Nichols algebra
over $\dm$, see Theorem \ref{teo:all-nichols-dm}.

\subsection{Nichols algebras of irreducible Yetter-Drinfeld
modules}\label{subsec:nichols-irred-dm}

In \cite[Table 2]{AF-alt-die}, it was determined the
dimension of Nichols algebras of
some irreducible Yetter-Drinfeld
modules over $\mathbb D_m$, with $m$
even. Here we complete the study in the
case $m = 4t \geq 12$, determining the
dimension of the Nichols algebras of
the irreducible Yetter-Drinfeld
modules coming from
the remaining two conjugacy classes $\Oc_{x}$ and $\oc_{xy}$.

\subsubsection{The conjugacy class of
$y^{n}$.}\label{subsubsec:y^n} Since
$y^n$ is central, the conjugacy class and the centralizer of $y^n$
in $\dm$ are $\oc_{y^n}=\{y^n\}$ and $C_{\dm}(y^n)=\dm$,
respectively. The irreducible representations of $\dm$ are
well-known and they are of degree 1 or 2. Explicitly, there are:
\begin{itemize}
  \item[(i)] $n-1=\frac{m-2}{2}$ irreducible representations of degree 2.
 Set $\omega$ an $m$-th primitive root of 1;
they are given by $\rho_\ell:\dm\to\GL(2,\k)$,
\begin{align}\label{eq:repdeg2}
\rho_\ell (x^ay^b)=\begin{pmatrix} 0 & 1 \\
                  1 & 0 \\
                \end{pmatrix}^a
 \begin{pmatrix}
                  \omega^\ell & 0 \\
                  0 & \omega^{-\ell} \\
                \end{pmatrix}^b, \quad \text{ $1\leq \ell<n$}.
\end{align}
\item[(ii)] 4 irreducible representations of degree 1.
They are given by Table \ref{tabladmpar}.
\end{itemize}
\begin{table}[h]
\begin{center}
\begin{tabular}{|c|c|c|c|c|c|}
\hline $\sigma$& 1 & $y^n$ & $y^b$, $1\leq b\leq n-1$ & $x$ & $xy$

\\ \hline  \hline $\chi_1$ &1&1&1&1&1

\\ \hline $\chi_2$ &1&1&1&$-1$&$-1$

\\ \hline $\chi_3$ &1&$(-1)^n$& $(-1)^b$ &1&$-1$

\\ \hline $\chi_4$ &1&$(-1)^n$& $(-1)^b$ &$-1$&1\\

\hline
\end{tabular}
\end{center}
\caption{One-dimensional irreducible
representations of $\dm$, with $m=2n$ even.}\label{tabladmpar}
\end{table}

Let $\rho$ be an irreducible representation of $C_{\dm}(y^n)=\dm$.
Since $n$ is even, if $\deg\rho=1$ or $\rho=\rho_\ell$ as in
\eqref{eq:repdeg2} with $\ell$ even, then $\dim\toba(\oc_{y^n},
\rho)=\infty$. Indeed, here $q_{y^ny^n}\neq -1$ and
Lemma \ref{le:odd} applies. In the cases when $\rho=\rho_\ell$ as in
\eqref{eq:repdeg2} with $\ell$ odd, we have that $q_{y^ny^n}= -1$.
Let $M_{\ell} = M(\oc_{y^n},\rho_\ell)$,
then we have the following lemma.

\begin{lema}\label{lema:ml} \cite[Thm. 3.1 (b) (i)]{AF-alt-die}
$\toba(\oc_{y^n}, \rho_{\ell}) \simeq \bigwedge
M_{\ell}$, for all $\ell$ odd with $1\leq \ell < n$. In particular,
 $\dim \toba(\oc_{y^n},
\rho_{\ell})= 4$. \qed
\end{lema}
Notice that there are $t$ irreducible Yetter-Drinfeld modules with
support $\oc_{y^n}$ such that its Nichols algebras is finite-dimensional.

\subsubsection{The conjugacy class of $y^{i}$, $1\leq i\leq n-1$.}
\label{subsubsec:y^j} The conjugacy class and the centralizer of
$y^i$ in $\dm$ are $\oc_{y^i}=\{y^i,y^{-i}\}$ and
$C_{\dm}(y^i)=\langle y \rangle \simeq \Z/m$, respectively. The
group of characters of $C_{\dm}(y^i)$ is
$$
\widehat{\quad C_{\dm}(y^i) = \, }  \{\chi_{(k)} \, | \, 1\leq k \leq m-1\},
$$
where $\chi_{(k)}(y):=\omega^k$, with
$\omega$ an $m$-th primitive root of 1. Let
$M_{i,k} = M(y^{i},\chi_{(k)})$.
Since $\oc_{y^i}$ is real, if $\chi_{(k)}(y^i)\neq -1$,
then $\dim\toba(\oc_{y^i},\chi_{(k)})=\infty$, by Lemma \ref{le:odd}.
Assume that $\chi_{(k)}(y^i)=-1$; this amounts to:
there exists $r$, with $r$ odd and $1\leq r\leq m-3$, such that
$ik= r n$.

Let
$\mathcal N_i=  \{k \, | \, 0\leq k \leq m-1, \chi_{(k)}(y^i)=-1 \}$.
Then, for every $i$ with $1\leq i \leq n-1$ there are
$\card\mathcal N_i$ irreducible Yetter-Drinfeld modules
with support $\oc_{y^i}$ and $\dim \toba(\oc_{y^i}, \chi_{(k)}) <\infty$.
Let $\omega \in \k$
be a primitive $m$-th root of 1. We define
$J=\{(i,k):\ \omega^{ik}=-1, 1\leq i\leq n-1, 1\leq k\leq m-1\}$.
\begin{rmk}
Notice that if $(i,m)=1$, then $\mathcal N_i=\{n\}$. Also,
\begin{itemize}
  \item if $i=2$, then $\mathcal N_2=\{t, 3t\}$;
  \item if $i=3$, then $\mathcal N_3=\{n\}$ if $3\not |\,t$, whereas
$\mathcal N_3=\{2u,6u,10u\}$ if $t=3u$;
  \item if $i=4$, then
 $\mathcal N_4=\emptyset$ if $2\not |\,t$, whereas $\mathcal
N_4=\{u, 3u,5u,7u\}$ if $t=2u$.
\end{itemize}
\end{rmk}


\begin{lema}\label{lema:mik} \cite[Thm. 3.1 (b) (ii)]{AF-alt-die}
$\toba(\oc_{y^i},  \chi_{(k)}) \simeq \bigwedge
M_{i,k}$, for all $(i,k)\in J$. In particular,
 $\dim \toba(\oc_{y^i},
\chi_{(k)})= 4$. \qed
\end{lema}

\subsubsection{The conjugacy classes of $x$ and $xy$.}
\label{subsubsec:x:xy}
We show that these two conjugacy
classes
give rise to infinite-dimensional Nichols algebras.

\begin{lema}
The classes $\oc_x$ and $\oc_{xy}$ are of type $D$. Hence
$\dim \toba(\oc_{x},\rho)$ and
$\dim \toba(\oc_{xy},\eta)$ are infinite for all
$\rho \in \widehat{C_{\dm}(x)}$ and $\eta \in
\widehat{C_{\dm}(xy)} $.
\end{lema}

\pf Since the classes $\oc_x$ and $\oc_{xy}$ are isomorphic as
racks to the dihedral rack $\mathcal D_{n}$, the result follows
from Lemma \ref{lema:d2m}. \epf

\begin{table}[ht]
\begin{center}
\begin{tabular}{|p{3.8cm}|p{2,2cm}|p{3cm}|p{1,7cm}|}
\hline {\bf Conjugacy class} &    {\bf Centralizer} & {\bf Rep.} &
$\dim \toba(V)$

\\ \hline  $e$   & $\dm$ & any & $\infty$
\\ \hline  $\Oc_{y^{\frac m2}}= \{y^{\frac m2}\}$, \newline
 $ \mid \Oc_{y^{\frac m2}} \mid = 1$   & $\dm$ &
$\chi_1$, $\chi_2$, $\chi_3$, $\chi_4$,
\newline $\rho_{\ell}$, $\ell$ even  & $\infty$
\\ \cline{3-4}     &  &  $\rho_\ell$, $\ell$ odd&   4
\\ \hline  $\Oc_{y^{i}}= \{y^{\pm i}\}$, $i\neq 0, \frac m2$,
\newline $ \mid \Oc_{y^{i}} \mid =2$   & $\Z/m \simeq \langle
y \rangle$ & $\chi_{(k)}$,  $\omega_m^{ik} = -1$ & 4
\\ \cline{3-4}     &  & $\chi_{(k)}$,
$\omega_m^{ik} \neq -1$ & $\infty$
\\ \hline  $\Oc_{x}= \{xy^{j}: j\text{ even}\}$
\newline $ \mid \Oc_{x} \mid = \frac{m}{2}$
& $\Z/2 \times \Z/2 \simeq$ \newline $\langle x\rangle \oplus
\langle y^{\frac m2} \rangle$ & $\varepsilon \otimes \varepsilon$,
$\varepsilon \otimes \sgn$,
\newline $\sgn \otimes \sgn$, $\sgn \otimes \varepsilon$ &$\infty$

\\ \hline  $\Oc_{xy}= \{xy^{j}:j\text{ odd}\}$
\newline $ \mid \Oc_{xy} \mid = \frac{m}{2}$
& $\Z/2 \times \Z/2 \simeq$ \newline $\langle xy \rangle \oplus
\langle y^{\frac m2} \rangle$ & $\varepsilon \otimes \varepsilon$,
$\varepsilon \otimes \sgn$,
\newline
$\sgn \otimes \sgn$, $\sgn \otimes \varepsilon$ &$\infty$

\\ \hline
\end{tabular}
\end{center}
\caption{$\mathbb D_{m}$, $m=4t$ with $t\geq 3$.}\label{tabladnpar:4t}
\end{table}

\subsection{Nichols algebras of arbitrary
Yetter-Drinfeld modules}
In this subsection we determine all finite-dimensional Nichols algebras in $\yddm$.
Specifically, we prove
that they are all exterior algebras over some
Yetter-Drinfeld modules. For such a module, we write
$\toba(M)$ instead of $\bigwedge M$ to enphasize the
Yetter-Drinfeld module structure.

\subsubsection{Nichols algebras over the family
$\{M_{i,k}\}$}\label{subsec:nichols-family-ik}
Recall $ M_{i,k} = M(\oc_{y^i}, \chi_{(k)})$,
with $1\leq i \leq n-1$, $0\leq k \leq m-1$.
We define an equivalence
relation in $J$, see Subsection \ref{subsubsec:y^j}, by
\begin{equation}\label{eq:rel-equiv-I}
(i,k) \sim (p,q)\mbox{ if }
\omega^{iq + pk}= 1.
\end{equation}
In such a case, one can prove that 
$\omega^{pk}=\omega^{iq}=-1$\footnote{Write $n = (i,n)h$. 
Since $n| ik$, one has that $h|k$. As
$n | iq + pk$, we have that $(i,n)|pk$ and thus $n| pk$. Then prove
that $pk \equiv n \mod m$.}. We
denote by $[i,k]=\{(p,q)\in J:\
(p,q) \sim (i,k)\}$ the class of $ (i,k)$ under this
equivalence.

\begin{prop}\label{prop:dimfin-irred-1}
Let
$M = M_{i_{1},k_{1}}\oplus \cdots \oplus M_{i_{r},k_{r}}$
with $(i_{s},k_{s}) \in J$ for all $1\leq s\leq r$. Then
$ \dim \toba(M) <\infty $ if and only if $(i_{p},k_{p}) \sim (i_{q},k_{q})$
for all $1\leq p,q\leq r$. In such a case, $\toba(M) \simeq
 \bigwedge M $ and
$ \dim \toba(M) =4^{r} $.
\end{prop}

\pf  Assume first $r=2$.
Let $ (i,k),(p,q) \in J $ and consider
$M_{i,k}$ and $M_{p,q}$. Then $\oc_{y^i}=\{\sigma_1:=y^i,
\sigma_2:=y^{-i}\}$, $\oc_{y^p}=\{\tau_1:=y^p,
\tau_2:=y^{-p}\}$ and
$\chi_{(k)}(y^i)=-1=\chi_{(q)}(y^p)
$.
Set $g_1=h_1=1$ and $g_2=h_2=x$, then
$$g_1 y^i g_1^{-1}=\sigma_1,\quad  g_2 y^i g_2^{-1} =
\sigma_2, \quad h_1 y^p h_1^{-1}=\tau_1, \quad h_2 y^p h_2^{-1} =\tau_2.$$

Consider now the Yetter-Drinfeld module $M=M_{i,k}\oplus M_{p,q}$.
As a vector space $M=\k$-span of $\{g_1,g_2,h_1,h_2\}$. The
braiding $c$ in $M$ is given by $c|_{M_{i,k}\oplus M_{i,k}}
=c_{M_{i,k}}$, $c|_{M_{p,q}\otimes M_{p,q}}=c_{M_{p,q}}$, and
$$
\begin{aligned}
c(g_1\ot h_1)&=\chi_{(q)}(y^i) \,h_1\ot g_1,
\quad  &c(g_1\ot h_2)&=\chi_{(q)}(y^{-i}) \,h_2\ot g_1,\\
c(g_2\ot h_1)&=\chi_{(q)}(y^{-i}) \,h_1\ot g_2,
\quad &c(g_2\ot h_2)&=\chi_{(q)}(y^{i}) \,h_2\ot g_2,\\
c(h_1\ot g_1)&=\chi_{(k)}(y^p) \,g_1\ot h_1,
\quad &c(h_1\ot g_2)&=\chi_{(k)}(y^{-p}) \,g_2\ot h_1,\\
c(h_2\ot g_1)&=\chi_{(k)}(y^{-p}) \,g_1\ot h_2,
\quad &c(h_2\ot g_2)&=\chi_{(k)}(y^{p}) \,g_2\ot h_2.
\end{aligned}
$$
Thus $M$ is a diagonal vector space whose matrix of coefficients is
$$\mathcal Q=\begin{pmatrix}
-1  & -1& \chi_{(q)}(y^i) & \chi_{(q)}(y^{-i}) \\
-1  & -1 & \chi_{(q)}(y^{-i}) & \chi_{(q)}(y^{i}) \\
 \chi_{(k)}(y^p)  & \chi_{(k)}(y^{-p}) & -1 & -1 \\
  \chi_{(k)}(y^{-p}) & \chi_{(k)}(y^p) & -1 & -1 \\
\end{pmatrix}.$$

Let
$\lambda:=\chi_{(q)}(y^i)  \chi_{(k)}(y^p) =
\omega^{iq+pk}$.
If $\lambda\neq 1$, then $(i,k)\nsim (p,q) $
and
$\dim\toba(M)=\infty$, by \cite{H2}, since
the generalized Dynkin diagram associated to $M$
is given by Figure \ref{fi:rombo}.
\begin{figure}[ht]
\vspace{1cm}
\begin{align*}
\setlength{\unitlength}{1.4cm}
\begin{picture}(1,0)
\put(0,0){\circle*{.15}} \put(1,1){\circle*{.15}}
\put(2,0){\circle*{.15}} \put(1,-1){\circle*{.15}}
\put(0,0){\line(1,1){1}} \put(0,0){\line(1,-1){1}}
\put(1,1){\line(1,-1){1}} \put(1,-1){\line(1,1){1}}
\put(0.2,.7){$\lambda$}
\put(1.6,.7){$\lambda^{-1}$}
\put(-0.05,-0.7){$\lambda^{-1}$}
\put(1.6,-0.7){$\lambda$}
\put(-0.4,-0.07){\small{-1}} \put(2.2,-0.07){\small{-1}}
\put(0.9,1.2){\small{-1}} \put(0.9,-1.3){\small{-1}}
\end{picture}\qquad \qquad
\end{align*}
\vspace{0.7cm}\caption{}\label{fi:rombo}
\end{figure}

If $\lambda= 1$, i.~e.  $\omega^{iq+pk}=1$ then $(p,q) \sim
(i,k)$. In such a case, $\omega^{pk}=\omega^{iq}=-1$, see
the paragraph after \eqref{eq:rel-equiv-I}, 
and whence $\toba(M)=\bigwedge M$, since
the braiding in $M$ is $c=-$flip; in particular $\dim (M) = 16$.

\smallbreak Assume $r\geq 2$ and let
$M = M_{i_{1},k_{1}}\oplus \cdots \oplus M_{i_{r},k_{r}}$
with $(i_{s},k_{s}) \in J$ for all $1\leq s\leq r$. In particular,
$\omega^{i_{s}k_{s}}=-1$ for all $1\leq s\leq r$. If
there exist $p$,
$q$, $1\leq p,q\leq r$ such that
    $(i_{p},k_{p}) \nsim (i_{q},k_{q})$, i.~e.
    $\chi_{(k_{q})}(y^{i_p})  \chi_{(k_{p})}(y^{i_{q}})
    = \omega^{i_{p}k_{q} + i_{q}k_{p}}
    \neq 1$, then $\dim \toba(M_{i_{p},k_{p}})
\oplus M_{i_{q},k_{q}})=\infty$ as
    above, which implies that $\dim \toba(M)=\infty$.
Thus $(i_{p},k_{p}) \sim (i_{q},k_{q})$ for all $1\leq p,q\leq r$
and $\chi_{(k_{q})}(y^{i_p}) \chi_{(k_{p})}(y^{i_{q}})=
\omega^{i_{p}k_{q} + i_{q}k_{p}}= 1$. As before,
$\omega^{i_{q}k_{p}}=\omega^{i_{p}k_{q}}=-1$, which implies that
$\toba(M)=\bigwedge M$, since the braiding in $M$ is $c=-$flip; in
particular $\dim\toba(M)=4^r$. \epf

\begin{definition}\label{def:I}
Let $$\II =\left\{ I=\coprod_{s=1}^{r}\{(i_{s},k_{s})\}: 
(i_{s},k_{s}) \in J \text{ and }(i_{s},k_{s})\sim (i_{p},k_{p}),
\ 1\leq s,p\leq r\right \}.$$ For $I\in \II$, we define $M_{I}=
\bigoplus_{(i,k)\in I}M_{i,k}$.
 \end{definition}
By Proposition \ref{prop:dimfin-irred-1}, we have
$\toba(M_{I})\simeq \bigwedge M_{I}$  and
$\dim\toba(M_{I}) = 4^{|I|}$.
\begin{obs}
Denote by $a_{i,k},
b_{i,k}$, $(i,k) \in I$ the
primitive elements that generate $\toba(M_{I})$.
Then, the Yetter-Drinfeld
module structure is given by
\begin{align}
\label{eq:prim:a:b:1}x\cdot a_{i,k} & =  b_{i,k},\quad & y\cdot a_{i,k}
& =  \omega^{k} a_{i,k},
\quad 
& \delta(a_{i,k}) & =   y^{i}\ot a_{i,k},\\
\label{eq:prim:a:b:2} x\cdot b_{i,k}  & =  a_{i,k},
\quad & y\cdot b_{i,k}  & = \omega^{-k} b_{i,k}
,\quad
& \delta(b_{i,k}) &= y^{-i}\ot b_{i,k}.
\end{align}
\end{obs}

\subsubsection{Nichols algebras over the family
$\{M_{\ell}\}$}\label{subsec:nichols-family-ell}
Recall that $M_{\ell} = M(\oc_{y^n}, \rho_{\ell})$.
In this subsection we study Nichols algebras over
sums of irreducible Yetter-Drinfeld modules
isomorphic to $M_{\ell}$,
with $1\leq \ell <n$, $\ell$ odd.

\begin{prop}\label{prop:dimfin-irred-2}
Let
$M = M_{\ell_{1}}\oplus \cdots \oplus M_{\ell_{r}}$
with $1\leq \ell_{s}<n$ odd numbers.
Then
$\toba(M) \simeq
\bigwedge M$ and
$ \dim \toba(M) =4^{r} $.
\end{prop}

\pf It suffices to show the braiding $c$ in $M$ is $c=-$flip. Let
$1\leq p,q\leq r $ and denote by $v_{1}, v_{2}$ and $w_{1},w_{2}$
the linear generators of $M_{\ell_{p}}$ and $M_{\ell_{q}}$,
respectively. Then $c = -$flip in $M_{\ell_{p}}\oplus
M_{\ell_{q}}$. Indeed, we know that $c|_{M_{\ell_{p}}\otimes
M_{\ell_{p}}} = -$flip and $c|_{M_{\ell_{q}}\otimes M_{\ell_{q}}}
= -$flip, by Lemma \ref{lema:ml}, and
\begin{align*}
c(v_{1}\otimes w_{1}) & =  y^{n}\cdot w_{1}\otimes v_{1} =
\omega^{n\ell_{q}} w_{1}\otimes v_{1}
= (-1)^{\ell_{q}}w_{1}\otimes v_{1} =
-w_{1}\otimes v_{1},\\
c(v_{1}\otimes w_{2}) & =  y^{n}\cdot w_{2}\otimes v_{1} =
\omega^{-n\ell_{q}} w_{2}\otimes v_{1}
= (-1)^{-\ell_{q}}w_{2}\otimes v_{1} =
-w_{2}\otimes v_{1},\\
c(v_{2}\otimes w_{1}) & =  y^{n}\cdot w_{1}\otimes v_{2} =
\omega^{n\ell_{q}} w_{1}\otimes v_{2}
= (-1)^{\ell_{q}}w_{1}\otimes v_{2} =
-w_{1}\otimes v_{2},\\
c(v_{2}\otimes w_{2}) & =  y^{n}\cdot w_{2}\otimes v_{2} =
\omega^{-n\ell_{q}} w_{2}\otimes v_{2}
= (-1)^{-\ell_{q}}w_{2}\otimes v_{2} =
-w_{2}\otimes v_{2},
\end{align*}
by straightforward computations.\epf

\begin{definition}\label{def:L}
Let $$\Lc = \left\{L=\coprod_{s=1}^{r} \{\ell_{s}\}: 
1\leq \ell_{1},\ldots, \ell_{r}<n\text{ odd numbers}\right\}$$ 
For $L\in \Lc$, we define $M_{L} = \bigoplus_{\ell \in L} M_{\ell}$.
\end{definition}

By Proposition \ref{prop:dimfin-irred-2}, we have:
$\toba(M_{L})\simeq \bigwedge M_{L} $ and $\dim
\toba(M_{L})=4^{|L|}$.

\begin{obs}
Denote by $c_{\ell}, d_{\ell}$ with $\ell \in L$ the primitive
elements that generate $\toba(M_{L})$.
Then, the  Yetter-Drinfeld module structure
is given by
\begin{align}
\label{eq:prim:c:d:1}x\cdot c_{\ell}  &= d_{\ell},  \quad & y\cdot c_{\ell}  
& = \omega^{\ell} c_{\ell},
\quad & \delta(c_{\ell}) &  =   y^{n}\ot c_{\ell},\\
\label{eq:prim:c:d:2}  x\cdot d_{\ell}  & =  c_{\ell}, 
\quad & y\cdot d_{\ell}  & =  \omega^{-\ell} d_{\ell},
\quad & \delta(d_{\ell}) & = y^{n}\ot d_{\ell}.
\end{align}
\end{obs}

\subsubsection{Nichols algebras over mixed
families}\label{subsec:nichols-family-mixed}

\begin{prop}\label{prop:dimfin-irred-3}
Let $M_{i,k,\ell} = M_{i,k}\oplus M_{\ell}$
with $(i,k) \in J$ and $1\leq \ell < n$ be an odd number.
Then
$\dim \toba(M_{i,k,\ell}) <\infty $ if and only if
$k$ is odd and $(i,\ell) \in J$.
In such a case,
$ \toba(M_{i,k,\ell})\simeq \bigwedge M_{i,k,\ell}$
and $ \dim \toba(M_{i,k,\ell}) =16 $.
\end{prop}

\pf Let $c:M_{i,k,\ell}\otimes M_{i,k,\ell}\to M_{i,k,\ell}\otimes
M_{i,k,\ell}$ be the braiding of $ M_{i,k,\ell} $. As before, it
suffices to show that $c= -$flip. Denote by $g_{1}= 1, g_{2}=x$
and $v_{1},v_{2}$ the linear generators of $M_{i,k}$ and
$M_{\ell}$, respectively. Then by Lemmata \ref{lema:mik} and
\ref{lema:ml}, we have that $c|_{M_{i,k}\otimes M_{i,k}} = -$flip
and $c|_{M_{\ell}\otimes M_{\ell}} = -$flip. Thus $c$ is
determined by the values
\begin{align*}
c(g_{1}\otimes v_{1}) & =  y^{i}\cdot v_{1}\otimes g_{1} =
\omega^{i\ell} v_{1}\otimes g_{1},\\
c(v_{1}\otimes g_{1}) & =  y^{n}\cdot g_{1}\otimes v_{1} =
\chi_{(k)}(y^{n})g_{1}\otimes v_{1} = \omega^{nk}g_{1}\otimes v_{1} =
(-1)^{k}g_{1}\otimes v_{1},\\
c(g_{2}\otimes v_{1}) & =  y^{-i}\cdot v_{1}\otimes g_{2} =
\omega^{-i\ell} v_{1}\otimes g_{2}, \\
c(v_{1}\otimes g_{2}) & =  y^{n}\cdot g_{2}\otimes v_{1} =
\chi_{(k)}(y^{n})  g_{2}\otimes v_{1}  = \omega^{nk}g_{2}\otimes v_{1}  =
(-1)^{k}g_{2}\otimes v_{1},\\
c(g_{1}\otimes v_{2}) & =  y^{i}\cdot v_{2}\otimes g_{1} =
\omega^{-i\ell} v_{2}\otimes g_{1},\\
c(v_{2}\otimes g_{1}) & =  y^{n}\cdot g_{1}\otimes v_{2} =
\chi_{(k)}(y^{n})g_{1}\otimes v_{2} = \omega^{nk}g_{1}\otimes v_{2} =
(-1)^{k}g_{1}\otimes v_{2},\\
c(g_{2}\otimes v_{2}) & =  y^{-i}\cdot v_{2}\otimes g_{2} =
\omega^{i\ell} v_{2}\otimes g_{2}, \\
c(v_{2}\otimes g_{2}) & =  y^{n}\cdot g_{2}\otimes v_{2} =
\chi_{(k)}(y^{n})  g_{2}\otimes v_{2}  = \omega^{nk}g_{2}\otimes v_{2}
=(-1)^{k}g_{2}\otimes v_{2}.
\end{align*}
This implies that $M$ is a diagonal vector space with matrix of coefficients
$$\mathcal Q=\begin{pmatrix}
-1  & -1&  \omega^{i\ell} & \omega^{-i\ell} \\
-1  & -1 & \omega^{-i\ell}& \omega^{i\ell} \\
(-1)^{k}  & (-1)^{k} & -1 & -1 \\
 (-1)^{k}& (-1)^{k} & -1 & -1 \\
\end{pmatrix}.$$
Let  $\lambda:= (-1)^{k}\omega^{i\ell}$.
If $\lambda\neq 1$, then
  the generalized
  Dynkin diagram associated to $M_{i,k,\ell}$
  is given by Figure \ref{fi:rombo} and whence
  $\dim\toba(M)=\infty$, by \cite{H2}.
Therefore, in order to have $ \dim\toba(M_{i,k,\ell})<\infty $ we
must have that $\lambda = 1$, that is, $(-1)^{k} =
\omega^{i\ell}$. By assumption $\omega^{ik}=-1$, thus
$\omega^{ik\ell} = (-1)^{\ell} = -1$, because $\ell$ is
odd. But $-1= (\omega^{i\ell})^{k} = ((-1)^{k})^{k} =
(-1)^{k^{2}}$, thus $k$ must be odd and $\omega^{i\ell}= -1$,
\textit{i.e.} $ (i,\ell) \in J $. In such a case, the braiding in
$M_{i,k,\ell}$ is $c= -$flip and then
$\toba(M_{i,k,\ell})\simeq \bigwedge M_{i,k,\ell}$ and $\dim
(M_{i,k,\ell}) = 16$. \epf

For $I\in \II$ and $L\in \Lc$, 
define
$M_{I,L} =
\left(\bigoplus_{(i,k) \in I} M_{i,k}\right)\oplus
\left(\bigoplus_{\ell \in L} M_{\ell}\right) $.
The next result generalizes Proposition \ref{prop:dimfin-irred-3}
for arbitrary finite sums.

\begin{prop}\label{prop:nichols-mixed}
Let $I \in \II, L\in \Lc$ and assume that $k$ is odd for all $(i,k)\in
I$. Then $\dim \toba(M_{I,L}) <\infty$ if and only if $(i,\ell)
\in J$ for all $(i,k) \in I$, $\ell\in L$. In such a case,
$\toba(M_{I,L}) \simeq \bigwedge M_{I,L}$ and $\dim \toba(M_{I,L})
= 4^{\vert I \vert + |L|} $.
\end{prop}

\pf Denote by $a_{i,k}, b_{i,k}$ and $c_{\ell},d_{\ell}$ the
linear generators of $M_{i,k}$ and $M_{\ell}$, respectively, for
all $(i,k) \in I$, $\ell\in L$. Then by the proof of Propositions
\ref{prop:dimfin-irred-1}, \ref{prop:dimfin-irred-2} and
\ref{prop:dimfin-irred-3}, it follows that $ \dim \toba(M_{I,L}) $
is finite if and only $k$ is odd for all $(i,k) \in I$ and
$(i,\ell) \in J$ for all $(i,k) \in I$, $\ell \in L$. In such a
case, the braiding in $M_{I,L}$ is given by $-$flip and whence
$\toba(M_{I,L}) \simeq \bigwedge M_{I,L}$. \epf

\begin{obs}
Denote by $a_{i,k}, b_{i,k}, c_{\ell},d_{\ell}$ with $(i,k) \in I$
and $\ell \in L$ the primitive elements that generate
$\toba(M_{I,L})$.
Then, the Yetter-Drinfeld module structure is determined by
\eqref{eq:prim:a:b:1}, \eqref{eq:prim:a:b:2},
\eqref{eq:prim:c:d:1} and \eqref{eq:prim:c:d:2}.
\end{obs}

\begin{definition}\label{def:I-L}
We define
$$\K = \left\{ (I,L):\ I\in \II, L\in \Lc\text{ and }  
 k\text{ odd}, (i,\ell)
\in J\text{ for all }(i,k) \in I, \ell\in L\right\}. $$
\end{definition}
By Proposition \ref{prop:nichols-mixed}, for all
$(I,L) \in \K$, we have 
$\toba(M_{I,L}) \simeq \bigwedge M_{I,L}$ and $\dim \toba(M_{I,L})
= 4^{\vert I \vert + |L|} $.
\subsection{Proof of Theorem A}
Let $\toba(M)$ be a finite-dimensional Nichols algebra in $\yddm$.
Since $\yddm$ is semisimple, $M$ must be a finite direct sum of
irreducible Yetter-Drinfeld modules. Then the result follows from
Lemmata \ref{lema:mik}, \ref{lema:ml} and Propositions
\ref{prop:dimfin-irred-1}, \ref{prop:dimfin-irred-2},
\ref{prop:dimfin-irred-3}. Clearly, Nichols algebras over distinct
families are pairwise non-isomorphic, since they are generated by
the set of primitive elements which are non-isomorphic as
Yetter-Drinfeld modules. \qed

\section{Liftings of Nichols algebras over
dihedral groups}\label{sec:pointed-ha-dm}
In this section we describe
all finite-dimensional pointed
Hopf algebras over dihedral groups
$\DD_{m}$,
assuming that $m=4t$, $t\geq 3$.

Let $H$ be a Hopf algebra with bijective antipode and let
$B \in \ydh$ be a braided Hopf algebra.
From
$B$ and $H$ one can construct a new Hopf algebra
$B\#H$, called the Majid-Radford product or bosonization,
whose underlying vector space is $B\ot H$ and
the Hopf algebra structure is given by
\begin{align*}
(a\# h)(b\# k) & =   a(h_{(1)}\cdot b )\# h_{(2)} k,\\
\Delta(a\#h) &= a^{(1)}\# (a^{(2)})_{(-1)} h_{(1)}
\ot (a^{(2)})_{(0)}\# h_{(2)},
\end{align*}
for all $(a\# h), (b\# k) \in B\#H$,
where $\Delta_{B}(a) = a^{(1)}\ot a^{(2)}$
is the braided coproduct and $\delta_{B}(a) = a_{(-1)}\ot a_{(0)}$
is the coaction of $H$ on $B$.

\begin{obs}\label{obs:rel-liftings}
 Assume that $A$ is a
finite-dimensional pointed Hopf algebra with $A_{0} = \k G(A)$ and
let $\gr A=\bigoplus_{i>0} A_i/A_{i-1}$. Then $\gr A$ is a Hopf
algebra which is isomorphic to the bosonization $R \#\k G(A)$,
where $R = A^{\co \pi}$. If $a \in R(1)$ is a homogeneos primitive
element, \textit{i.e.} $\delta(a) = g\ot a$, $g \in G(A)$, then
$a\# 1 \in R\# \k G(A)$ is $(g,1)$-primitive. Indeed, by the
formula above we have $ \com (a\#1) = a^{(1)}\# (a^{(2)})_{(-1)}
\ot (a^{(2)})_{(0)}\# 1 = a\# 1 \ot 1\# 1 + 1\#g \ot a \# 1 $.
Consider now the projection, $\pi:A_1\to A_1/A_0$, which is in
particular a proyection of Hopf $\k G(A)$-bimodules, and denote by
$\sigma$ a section of Hopf $\k G(A)$-bimodules. Since $A_{1}/A_{0}
= A_{0}\oplus P(R)\# \k G(A)$, by \cite[Lemma 2.4]{AS0}, we have
that $a\#1 \in A_{1}/A_{0}$, and $\sigma(a\#1)$ is
$(g,1)$-primitive in $A$.
\end{obs}

The following is a key step for the classification, see
\cite[Prop. 5.4]{AS}, \cite[Thm. 7.6]{AS2}, 
\cite[Thm. 2.1]{AG2}, \cite[Thm. 3.1]{GG}, \cite[Thm. 2.6]{AGI}.
It agrees with a well-known conjecture \cite[Conj. 1.4]{AS}.

\begin{theorem}\label{teo:item-c}
Let $A$ be a finite-dimensional pointed Hopf algebra with
$G(A) = \DD_{m}$. Then $A$ is generated by
group-likes and skew-primitives.
\end{theorem}

\pf Since $\gr A = R\#\k \DD_{m}$, with $R= \bigoplus_{n\geq
0}R(n)$ the diagram of $A$, it is enough to prove that $R$ is a
Nichols algebra, since in such a case, $A$ would be generated by
$G(A)$ and skew-primitive elements. Let $S$ be the graded dual of $
R $. By \cite[Lemma 5.5]{AS}, $S$ is generated by $V=S(1)$ and $R$
is a Nichols algebra if and only if $P(S)=S(1)$, that is, if $S$
is itself a Nichols algebra.

Consider $\toba(V) \in
\yddm$. Since $V=R(1)^{*}= P(R)^{*}$ and
$\toba(P(R))$ is finite-dimensional, we have
by \cite[Prop. 3.2.30]{AG1} that $\toba(V)$ is also finite-dimensional and
by Theorem \ref{teo:all-nichols-dm},
$\toba(V)$ is isomorphic to an exterior algebra
$ \toba(M_{I}) $, $ \toba(M_{L})$ or $\toba(M_{I,L}) $,
with $I\in \II$, $L\in \Lc$ and $(I,L) \in \K$, respectively.
Moreover,
a direct computation shows that the elements $r$ in $S$
representing the quadratic relations are primitive and since 
the braiding is $-$flip, they
satisfy that $c(r\ot r) =r\ot r$. As $\dim S <\infty$,
it must be $r=0$ in $S$ and hence there
exists a proyective algebra map $\toba(V)
\twoheadrightarrow S$, which implies that
$S$ is generated by $S(1)$.
\epf

\subsection{Some liftings and quadratic relations}
We begin this subsection with
the following 
proposition 
that shows how to deform the relations
in the Nichols algebras to get a lifting.

Let $A$ be a finite-dimensional pointed Hopf algebra 
over $\k \DD_{m}$. Then by Theorem
\ref{teo:item-c}, we have that $\gr A = \toba(M)\# \k \dm$ and
its infinitesimal braiding $M$ is
isomorphic either to $M_{I}$ with $I \in \II$ and
$|I|>1$  or $I = \{(i,n)\}$, or $M_{L}$ with $L\in \Lc$, or
$M_{I,L}$ with $(I,L) \in \K$ and $|L|>0$, $|I|> 0$, see Subsections
\ref{subsec:nichols-family-ik}, \ref{subsec:nichols-family-ell}
and \ref{subsec:nichols-family-mixed}. 

From now on, we denote by $g,\ h$ the generators of $G(A)\simeq \DD_{m}$ with
$g^{2} = 1 = h^{m}$ and $ghg = h^{-1}$. 

For all $1\leq r,s<m$, let $M^{s}_{r} = 
\{a\in M:\ \delta(a) = h^{s}\ot a, h\cdot a = \omega^{r}a\}$.
Then $M = \bigoplus_{r,s} M^{s}_{r}$.
Following Remark
\ref{obs:rel-liftings}, we write  
$x = \sigma(a\#1)$
for the element in $A$ defined
by the Hopf $\k G(A)$-bimodule section $\sigma$.
In particular, 
if $a \in M_{r}^{s}$, then 
$x$ is $(h^{s},1)$-primitive and $h x h^{-1} = \omega^{r} x$.

\begin{prop}\label{prop:def-rel-lifting}
Let $A$ be a finite-dimensional pointed Hopf algebra 
with $G(A) = \dm$ and infinitesimal braiding $M$.
Let $a \in M^{s}_{r}, b \in M^{v}_{u}$ with
$1\leq r,s,u,v <m$ and 
denote $x= \sigma(a\#1)$, $y=\sigma(b\#1)$. Then
there exists $\lambda \in \k$ such 
that 
\begin{equation} \label{eq:def-rel-lifting}
xy + yx = \delta_{u,m-r}\lambda (1 - h^{s+v}). 
\end{equation}
In particular, if $x=y$ we have that $x^{2}= \delta_{r,n}\lambda' (1-h^{2s})$
with $\lambda' = \frac{\lambda}{2}$. 
\end{prop}

\pf
By Theorem \ref{teo:item-c}, $M$ is isomorphic either to
to $M_{I}$ with $I \in \II$ and
$|I|>1$  or $I = \{(i,n)\}$, or $M_{L}$ with $L\in \Lc$, or
$M_{I,L}$ with $(I,L) \in \K$ and $|L|>0$, $|I|> 0$. As  
$a \in M^{s}_{r}, b \in M^{v}_{u}$,  Propositions
\ref{prop:dimfin-irred-1}, \ref{prop:dimfin-irred-2} and \ref{prop:nichols-mixed} yield
that $\omega^{sr}= -1 = \omega^{uv}$ and $(r,s)\sim (u,v)$,
\textit{i.e.} $\omega^{rv + su} = 1$ and $\omega^{rv} = -1 =
\omega^{su}$. 

A straightforward computation yields that
the element $\alpha=
xy +yx$ is
$(h^{s+v},1)$-primitive. Indeed, 
\begin{align*}
& \com(\alpha)
=\com(xy +yx) \\
&\,=(x\ot 1+ h^{s}\ot  x)(y\ot 1 + h^{v}\ot y)+
(y\ot 1 + h^{v}\ot y)
(x\ot 1+ h^{s}\ot  x)\\
&\, = xy\ot 1 + xh^{v}\ot  y +
h^{s}y\ot x +
h^{s+v}\ot xy +yx\ot 1 +\\
& \quad \qquad+yh^{s}\ot x
+h^{v}x\ot y
+h^{s+v}\ot yx \\
&\,  = (xy+yx)\ot 1 +
h^{s+v}\ot (xy+yx)
 +(xh^{v}+h^{v}x)\ot y
 + (h^{s}y +yh^{s})\ot x\\
& \, = (xy+yx)\ot 1 + h^{s+v}\ot
(xy+yx).
 \end{align*}
If $s+v \equiv 0 \mod m$, then $\alpha$ is primitive.
Since $A$ is finite-dimensional, we must have that 
$\alpha = 0$. Suppose $s+v\not\equiv 0\mod m$.
Then, by Theorem \ref{teo:item-c}
there exist $(h^{s+v},1)$-primitive elements $x_{i,j} \in M^{i}_{j}$
with $i=s+v$ and $\lambda, \beta_{i,j} \in \k$ 
such that
$$  \alpha=
\lambda(1 - h^{s+v}) + \sum_{i=s+v, j}
\beta_{i,j}x_{i,j}.
$$
Conjugating on both sides by $h$ gives
$$
h\alpha h^{-1}   =
\omega^{r+u}\alpha 
= \lambda_{p,q,i,k}(1 - h^{s+v})
+ \sum_{i=s+v,j} \beta_{i,j}\omega^{j}x_{i,j},$$
which implies that $ \lambda=
\lambda\omega^{r+u}$ and $\beta_{i,j}\omega^{j}
=\beta_{i,j}\omega^{r+u}$ for all
$i,j$. If $r+u \not\equiv 0
\mod m$, then $ \lambda_{p,q,i,k}=0 $. Thus, to end the 
proof we need to show that
necessarily $\beta_{i,j} = 0 $ for all $i,j$. Suppose on the contrary 
that
$\beta_{i,j} \neq 0$ for some $i,j$. Then $j\equiv r+u\mod m$. In such a
case, as $i = s+v$ we have
$$-1=\omega^{ij}=
 \omega^{(s+v)(r+u)} = \omega^{sr + vu}
\omega^{su+vr}=1,$$
a contradiction. In conclusion, we must have
that $\alpha=
\delta_{u,m-r} \lambda(1 - h^{s+v})$.
\epf

As a direct consequence of Proposition \ref{prop:def-rel-lifting}
we get the following corollaries.
The first one shows that all
pointed Hopf algebras over $\dm$ whose
infinitesimal braiding is isomorphic to
$M_{I}$, with $I = \{(i,k)\} \in \II$, or $M_{L}$ 
with $L\in \Lc$, as in Subsections
\ref{subsec:nichols-family-ik} and
\ref{subsec:nichols-family-ell}, respectively,
are isomorphic to bosonizations.

\begin{cor}\label{prop:sin-lifting-ik-1}\label{prop:sin-lifting-m-l}
Let $A$ be a finite-dimensional pointed Hopf algebra with
 $G(A)\!=\DD_{m}$, such that its infinitesimal braiding $M$ is
isomorphic to 
$M_{I}$ with $I = \{(i,k)\} \in \II$ and $k\neq n$, or
 $ M_{L} $ with $L \in \Lc$. 
Then
$A \simeq \gr A\simeq \toba (M)\# \k \DD_{m}$.
\end{cor}

\pf
Suppose $M\simeq M_{I}$, with $I = \{(i,k)\} \in \II$ and $k\neq n$ and
denote $x = \sigma(a_{i,k}\#1)$, $y = \sigma(b_{i,k}\#1)$.
Then Proposition \ref{prop:def-rel-lifting} gives
$x^{2} = 0 = y^{2}$ and $xy+yx = \delta_{i,m-i} \lambda (1-h^{i+m-i}) = 0$
for any $\lambda \in \k^{\times}$. Thus $A \simeq \gr A\simeq 
\toba (M_{I})\# \k \DD_{m}$. Assume now that 
$ M_{L} $ with $L \in \Lc$ and 
denote $x = \sigma(c_{\ell}\#1)$, $y = \sigma(c_{\ell'}\#1)$
with $\ell, \ell' \in L$ and $e_{\ell}, e_{\ell'}$ in the
set of linear generators  $\{c_{\ell},d_{\ell}: \ell \in L\}$ of $M$.
As $x$ and $y$ are $(h^{n},1)$-primitive, a similar computation
as above shows that $x$, $y$ and $xy+yx$ are primitive. Hence 
 $A \simeq \gr A\simeq 
\toba (M_{L})\# \k \DD_{m}$ and the corollary is proved.
\epf

The following two corollaries give the explicit relations
that a lifting of a Nichols algebra over $\dm$ must satisfy.
 
\begin{cor}\label{cor:rel-lifting-ikI}
Let $A$ be a finite-dimensional pointed Hopf algebra with
$G(A)\!=\DD_{m}$, such that its infinitesimal braiding is
isomorphic to $M_{I} $, with $I\in \II$ and $|I| > 1$ or $I = \{(i,n)\}$.
Denote $x_{p,q}=\sigma(a_{p,q}\#1)$ and $y_{p,q}=\sigma(b_{p,q}\# 1)$
for all $(p,q)\in I$.
Then there exist two families of elements
in $\k$, $\lambda = (\lambda_{p,q,i,k})_{(p,q),(i,k) \in I}$,
and $\gamma = (\gamma_{p,q,i,k})_{(p,q),(i,k) \in I}$,
such that
 \begin{align}
\label{eq:1} x_{p,q}x_{i,k} + x_{i,k}x_{p,q} & =
\delta_{q,m-k} \lambda_{p,q,i,k}(1 - h^{p+i}),
\\ \label{eq:2}
 y_{p,q}y_{i,k} + y_{i,k}y_{p,q} & =
\delta_{q,m-k} \lambda_{p,q,i,k}(1 - h^{-p-i}),\\ \label{eq:3}
x_{p,q}y_{i,k} + y_{i,k}x_{p,q} &=
\delta_{q,k} \gamma_{p,q,i,k}(1 - h^{p-i}).
\end{align}\qed
\end{cor}

\begin{obs}
Note that the symmetry of relations \eqref{eq:1} and \eqref{eq:3} imply that
the families $\lambda= (\lambda_{p,q,i,k})_{(p,q),(i,k) \in I}$,
and $\gamma = (\gamma_{p,q,i,k})_{(p,q),(i,k) \in I}$
satisfy
\begin{equation}\label{eq:sat-par-1}
\lambda_{p,m-k,i,k}=  \lambda_{i,k,p,m-k}\quad
\text{and}\quad
\gamma_{p,k,i,k}= \gamma_{i,k,p,k}.
\end{equation}
\end{obs}

\begin{cor}\label{cor:rel-lifting-ikI-L}
Let $A$ be a finite-dimensional pointed Hopf algebra with
$G(A)\!=\DD_{m}$ such that its infinitesimal braiding is
isomorphic to $M_{I,L}$,  with $(I,L)\in \K$.
Denote $x_{p,q}=\sigma(a_{p,q}\#1)$, 
$y_{p,q}=\sigma(b_{p,q}\# 1)$,
$z_{\ell}= \sigma(c_{\ell}\# 1)$ and
$w_{\ell}= \sigma(d_{\ell}\# 1)$
for all $(p,q)\in I$, $\ell \in L$.
Then there exist
four families of elements
in $\k$,
$\lambda = (\lambda_{p,q,i,k})_{(p,q),(i,k) \in I}$,
$\gamma = (\gamma_{p,q,i,k})_{(p,q),(i,k) \in I}$,
$\theta = (\theta_{p,q,\ell})_{(p,q) \in I, \ell \in L}$,
and $\mu = (\mu_{p,q,\ell})_{(p,q)\in I, \ell \in L}$,
such that the following relations in $A$ hold:
 \begin{align}\label{eq:1b}
 x_{p,q}^{2} = 0 = y_{p,q}^{2},&\qquad
z_{\ell}w_{\ell\,'}+w_{\ell\,'}z_{\ell}= 0,\qquad\\ \label{eq:2b}
z_{\ell}z_{\ell\,'}+z_{\ell\,'}z_{\ell}= 0,&\qquad
w_{\ell}w_{\ell\,'}+w_{\ell\,'}w_{\ell}= 0\\ \label{eq:3b}
 x_{p,q}x_{i,k} + x_{i,k}x_{p,q} & =
\delta_{q,m-k} \lambda_{p,q,i,k}(1 - h^{p+i}),\qquad\\ \label{eq:4b}
 y_{p,q}y_{i,k} +y_{i,k}y_{p,q} & =
\delta_{q,m-k} \lambda_{p,q,i,k}(1 - h^{-p-i}),\qquad\\ \label{eq:5b}
x_{p,q}y_{i,k} +y_{i,k}x_{p,q} & =
 \delta_{q,m-k}\gamma_{p,q,i,k}(1 - h^{p-i}),\qquad\\ \label{eq:6b}
 x_{p,q}z_{\ell} + z_{\ell}x_{p,q} & =
 \delta_{q,m-\ell}\,\theta_{p,q,\ell}(1 - h^{n+p}),\qquad\\ \label{eq:7b}
 y_{p,q}w_{\ell} + w_{\ell}y_{p,q} & =
  \delta_{q,m-\ell}\,\theta_{p,q,\ell}(1 - h^{n-p}),\qquad\\ \label{eq:8b}
 x_{p,q}w_{\ell} + w_{\ell}x_{p,q}& =
 \delta_{q,\ell}\,\mu_{p,q,\ell}(1 - h^{n+p}),\qquad\\ \label{eq:9b}
 y_{p,q}z_{\ell} + z_{\ell}y_{p,q} & =
 \delta_{q,\ell}\,\mu_{p,q,\ell}(1 - h^{n-p}).\qquad
\end{align}\qed
\end{cor}

\begin{obs}\label{rmk:eq2}
As before, \eqref{eq:3b} and \eqref{eq:5b} imply
the equalities in \eqref{eq:sat-par-1}.
\end{obs}

\subsection{Quadratic algebras}
In this section we introduce two families of pointed
Hopf algebras which are defined by quadratic relations.
They are constructed by deforming the
relations on two families of Nichols algebras
in $\yddm$. Moreover,
we show that they are liftings of bosonizations
of Nichols algebras and belong to the family of
Hopf algebras that characterize pointed Hopf algebras
over $\DD_{m}$, $m=4t$, $t\geq 3$.

\bigbreak \noindent {\bf The families of parameters}. 
Let $I\in \II$ and $L\in \Lc$
be as in Definitions \ref{def:I} and \ref{def:L}, respectively, and let
$\lambda = (\lambda_{p,q,i,k})_{(p,q),(i,k) \in I}$, $\gamma =
(\gamma_{p,q,i,k})_{(p,q),(i,k) \in I}$, $\theta =
(\theta_{p,q,\ell})_{(p,q) \in I, \ell \in L}$, and $\mu =
(\mu_{p,q,\ell})_{(p,q)\in I, \ell \in L}$ be families of elements
in $\k$, satisfying the following conditions:
\begin{equation}\label{eq:cond-par-1}
\lambda_{p,m-k,i,k}=  \lambda_{i,k,p,m-k}\quad\text{and}\quad
\gamma_{p,k,i,k}= \gamma_{i,k,p,k}.
\end{equation}
In particular, $\theta$ and $\mu$ are families of free
parameters in $\k$.

\begin{definition}\label{liftings1}
For $I \in \II$, denote by
$A_{I}(\lambda, \gamma)$ the
algebra generated by the elements $g,h, x_{p,q},y_{p,q}$
with $(p,q) \in I$ satisfying
the following relations:
\begin{align}
g^{2} = 1 = h^{m},\qquad \qquad ghg= h^{m-1},\qquad\ \qquad\\
gx_{p,q}  = y_{p,q} g,\qquad  hx_{p,q} = \omega^{q} x_{p,q} h,
\qquad hy_{p,q} = \omega^{-q}y_{p,q}h,\\
x_{p,q}x_{i,k} +x_{i,k}x_{p,q}  =
\delta_{q,m-k} \lambda_{p,q,i,k}(1 - h^{p+i}),
\qquad\label{eq:x-x-a-I}\\
x_{p,q}y_{i,k} +y_{i,k}x_{p,q} =
\delta_{q,k} \gamma_{p,q,i,k}(1 - h^{p-i}).\qquad\label{eq:x-y-a-I}
\end{align}
\end{definition}
It is a Hopf algebra with its structure determined by
\begin{align*}
\com(g) &= g\ot g, \qquad\qquad\qquad\qquad\quad \com(h) = h\ot h,\\
\com(x_{p,q}) & = x_{p,q}\ot 1 + h^{p}\ot x_{p,q},\qquad
\com(y_{p,q}) = y_{p,q}\ot 1 + h^{-p}\ot y_{p,q}.
\end{align*}
Since it is generated by group-likes and skew-primitives, it is
pointed by \cite[Lemma 5.5.1]{M}. We will call the pair $
(\lambda,\gamma) $ a \emph{lifting datum} for $\toba(M_{I})$. We
set $\gamma = 0$ if $|I|=1$.

\begin{exa} If $I=\{(i,k)\}$ with $k\neq n$,
then by Corollary
\ref{prop:sin-lifting-ik-1}, the Hopf algebra defined
above is the bosonization $\toba(M_{I})\#\k \DD_{m}$. If
$k=n$
we obtain
the Hopf algebra $A_{i,n}(\lambda)$ generated
by the elements $g, h, x, y$ satisfying
\begin{eqnarray}
\nonumber g^{2} = 1 = h^{m},&\qquad ghg= h^{m-1},\\
\nonumber gx  = y g,&\qquad hx = -x h,
\qquad hy = -yh,\\
\nonumber x^{2}= \lambda (1 - h^{2i}),
&\quad y^{2}=
 \lambda(1 - h^{-2i}),\quad
xy + yx = 0.
\end{eqnarray}
It is a finite-dimensional pointed Hopf algebra
with its  structure given by
$$\com(g) = g\ot g,\quad
\com(h) = h\ot h,\quad
\com(x) = x\ot 1 + h^{i}\ot x,\quad
\com(y) = y\ot 1 + h^{-i}\ot y.$$
\end{exa}

\begin{definition}\label{def:liftings2}
For $(I,L) \in \K$, denote by
$B_{I,L}(\lambda,\gamma,\theta,\mu)$ the
algebra generated by  $g,h, x_{p,q},y_{p,q},z_{\ell},w_{\ell}$
satisfying
the relations:
\begin{eqnarray}
g^{2} = 1 = h^{m},&\qquad ghg= h^{m-1},\\
gx_{p,q}  = y_{p,q} g,&\qquad hx_{p,q} = \omega^{q} x_{p,q} h,
\\
gz_{\ell}  = w_{\ell} g,&\qquad hz_{\ell} = \omega^{\ell} z_{\ell} h,
\\
x_{p,q}^{2} = 0 = y_{p,q}^{2}&\qquad
z_{\ell}w_{\ell'}+w_{\ell'}z_{\ell}= 0\qquad
z_{\ell}z_{\ell'}+z_{\ell'}z_{\ell}= 0\\
& x_{p,q}x_{i,k} +x_{i,k}x_{p,q} =
\delta_{q,m-k} \lambda_{p,q,i,k}(1 - h^{p+i}),\qquad\\
&x_{p,q}y_{i,k} +y_{i,k}x_{p,q} =
 \delta_{q,k}\gamma_{p,q,i,k}(1 - h^{p-i}),\qquad\\
 &x_{p,q}z_{\ell} + z_{\ell}x_{p,q} =
 \delta_{q,m-\ell}\theta_{p,q,\ell}(1 - h^{n+p}),\qquad\\
 &x_{p,q}w_{\ell} + w_{\ell}x_{p,q} =
 \delta_{q,\ell}\mu_{p,q,\ell}(1 - h^{n+p}).\qquad
\end{eqnarray}
\end{definition}
It is a Hopf algebra with its structure determined by
\begin{align*}
\com(g) &= g\ot g, & \com(h) &= h\ot h,\\
\com(x_{p,q}) & = x_{p,q}\ot 1 + h^{p}\ot x_{p,q},&
\com(y_{p,q}) & = y_{p,q}\ot 1 + h^{-p}\ot y_{p,q},\\
\com(z_{\ell}) & = z_{\ell}\ot 1 + h^{n}\ot z_{\ell},&
\com(w_{\ell})& = w_{\ell}\ot 1 + h^{n}\ot w_{\ell}.
\end{align*}
Since it is generated by group-likes and skew-primitive elements,
it is pointed by \cite[Lemma 5.5.1]{M}. We call the 4-tuple $
(\lambda,\gamma,\theta,\mu) $ a lifting datum for
$\toba(M_{I,L})$.

\begin{exa} If $I=\{(i,k)\}$ and $L=\{\ell\}$
 with $1\leq k,\ell<m$ odd numbers
and $m-\ell \neq k$,
then
the Hopf algebra defined
above is the bosonization $\toba(M_{I,L})\#\k \DD_{m}$. If
$k=m-\ell$
we obtain
the Hopf algebra $B_{I,L}(\theta,\mu)$ generated
by the elements $g, h, x, y,z,w$ satisfying the
relations
\begin{eqnarray}\nonumber
g^{2} = 1 = h^{m},&\qquad ghg= h^{m-1},\\\nonumber
gx  = y g,&\qquad hx = \omega^{k}x h,
\qquad
gz  = w g,  \qquad hz =  \omega^{-k}zh,
\\\nonumber
x^{2}= 0 = y^{2},&\qquad z^{2} = 0 = \omega^{2},\qquad
xy + yx = 0,\\\nonumber
zw +wz= 0,&\quad xz+zx= \theta (1 - h^{n+i}),\quad
xw+wx= \mu(1 - h^{n+i}).\nonumber
\end{eqnarray}
\end{exa}

As we have seen in Section \ref{sec:nichols-dihedral},
finite-dimensional Nichols algebras in $\yddm$ are exterior
algebras, see Theorem \ref{teo:all-nichols-dm}. Also, the Hopf
algebras $A_{I}(\lambda,\gamma) $ and $B_{I,L}(\lambda,\gamma,
\theta,\mu)$ defined above are quadratic algebras for 
all $I\in \II$ and $(I,L)\in \K$. Our next goal
is to show that if $H$ is a lifting of a finite-dimensional
Nichols algebra in $ \yddm $, then $H$ is isomorphic to a
quadratic algebra defined above for some lifting data, and
conversely, these Hopf algebras together with the bosonizations
are all liftings of finite-dimensional Nichols algebras in $ \yddm
$. First we work on quadratic algebras to obtain a bound on the
dimensions of $A_{I}(\lambda,\gamma)$ and $B_{I,L}(\lambda,\gamma,
\theta,\mu)$. We follow \cite{GG} for our exposition.

\par Let $ W $ be a finite-dimensional
vector space and let $T(W) = \oplus_{n\geq 0}W^{\ot n}$ be the graded
tensor algebra with the induced increasing filtration
$F^i:=\oplus_{j\leq i}W^{\ot j}$. Let $R\subset W\ot W$ be a
subspace and denote by $J(R)$ the two-sided ideal of $T(W)$ generated
by $R$. A (homogeneous) quadratic algebra $Q(W,R)$ is the quotient
$T(W)/J(R)$. Analogously, for a subspace $P\subset F^2=\k\oplus
W\oplus W\ot W$, we denote by $J(P)$ the two-sided ideal in $T(W)$
generated by $P$. A (nonhomogeneous) quadratic algebra $Q(W,P)$ is
the quotient $T(W)/J(P)$.

\par Let $A=Q(W,P)$ be a nonhomogeneous quadratic algebra. It inherits an
increasing filtration $A_n$ from $T(W)$ given by
$A_n=F^n/(J(P)\cap F^n)$ and the associated graded algebra is
 $\gr A=\oplus_{n\geq 0} A_n/A_{n-1}$, where $A_{-1}=0$.
Consider now the projection $\pi: F^2\to W\ot W$ with kernel $F^1$
and set $R=\pi(P)\subset W\ot W$. Let $B=Q(W,R)$ be the
homogeneous quadratic algebra defined by $R$. If $P\cap W = 0$,
then we have an epimorphism $\rho: B\to \gr A$. Indeed, let
$\rho': T(W)\to \gr A$ be the graded algebra map induced by
$W\hookrightarrow A_1 \twoheadrightarrow A_1/A_{0}$. Suppose $x\in
R\subset W^{\ot 2}$, then there exist $x_0\in\k$, $x_1\in W$ such
that $x-x_1-x_0\in P$ and therefore $ x=x_1+x_0\in F^2/( J(P)\cap
F^2)=A_2$. Thus $\rho'(x)=0\in A_2/A_1$, and whence $\rho'$
induces $\rho:B=T(W)/J(R)\twoheadrightarrow \gr A$.

\par Let $G$ be a finite group and suppose that
$\toba(V)$ is a finite-dimensional Nichols algebra in
$\ ^{\k G}_{\k G}\mathcal{YD}$
which is given by quadratic relations, \textit{i.e.}
$\toba(V) = \widehat{\toba}_{2}(V)$. Denote by
$p_{i}(x_{j_{1}},\ldots, x_{j_{k}}) = 0 $ these quadratic relations,
where $x_{j} \in P(V)$ are $g_{j}$-homogeneous
elements with $g_{j} \in G$ and $p_{i}, i\in I$, are
finitely many quadratic polynomials with coefficients in $ \k$.

\par Assume $H$ is a
Hopf algebra containing $\k G$ as a
Hopf subalgebra which is generated by $\k G$  and  $P(V)$
such that
$x_{j}$ is $(g_{j},1)$-primitive for all $x_{j}\in P(V)$,
$gx_{j}g^{-1} = g\cdot x_{j} $ for all $g\in G, x_{j}\in P(V)$
and $p_{i}(x_{j_{1}},\ldots, x_{j_{k}}) =
\lambda_{i,j_{1},\ldots, j_{k}}(1-g_{i,j_{1},\ldots, j_{k}})$ for some
$\lambda_{i,j_{1},\ldots, j_{k}} \in \k$ and
$g_{i,j_{1},\ldots, j_{k}} \in G$.
The next lemma is a slight generalization of \cite[Prop. 4.2]{GG}.

\begin{lema}\label{lema:quadratic-bound}
$\dim H \leq \dim \toba(V)\vert G\vert$.
\end{lema}

\pf $H$ is the nonhomogeneous quadratic algebra
$Q(W,P)$ defined by $W$ and $P$, for $W=\k\{x_j, H_g : x_{j}\in P(V),
g\in
G\}$ and $P\subset \k\oplus W\oplus W\ot W$ the subspace generated
by
\begin{align*}
\{H_e-1&, H_g\ot H_t-H_{gt}, H_g\ot x_j- g\cdot x_{j}\ot
H_g,\\ &\qquad   p_{i}(x_{j_{1}},\ldots, x_{j_{k}}) -
\lambda_{i,j_{1},\ldots, j_{k}}(1-g_{i,j_{1},\ldots, j_{k}})\}.
\end{align*}
Let $R=\pi(P)$.
Explicitly, $R\subset W\ot W$ is the subspace generated by
$
\{H_g\ot H_t, H_g\ot x_j- g\cdot x_{j}\ot
H_g, p_{i}(x_{j_{1}},\ldots, x_{j_{k}})\}
$.
Let $B=Q(W,R)$ be the homogeneous quadratic algebra defined by $W$
and $R$. Then $B\cong
\toba(V)\# Y_{G}$, where $Y_{G}$ is the algebra linearly
spanned by the set $\{1, y_g:g\in G\}$ with unit 1 and
multiplication table given by $ y_gy_t=0$ for all $g,t\in G$ and
$\#$ stands for the commutation relation
$(1\# y_{g})(x_{j}\# 1) = g\cdot x_{j}\# y_{g} $,
$ (1\# 1)(x_{j}\# 1) = x_{j}\# 1 $ for all
$g\in G$, $x_{j}\in P(V)$. Thus by the preceeding discussion,
there exists an epimorphism $\rho:
\toba(V)\#  Y_{G} \twoheadrightarrow \gr H$.

\par Since $P\cap F^1=\k\{H_e-1\}$,
by \cite[Lemma 4.1]{GG} we have $\rho(H_e)=0$ and whence there
exists an epimorphism $\rho_e:B/By_eB\to \gr H$. The commutation
relation and the fact that the elements $\{y_g\}_{g\in G}$ are
pairwise orthogonal, give $By_eB=\toba(V)y_e\subset B$. This
implies $\dim B^n -\dim (\toba(V)^n y_e)\geq \dim\gr H^n$ and
since $\dim B^n=\dim \toba(V)^n (|G|+1)$, we have $\dim \toba(V)^n
|G|\geq \dim \gr H^n$ and consequently $\dim H \leq
\dim\toba(V)\vert G\vert$. \epf

The next corollary follows inmediately.

\begin{cor}\label{cor:dim-bound} For all $I\in \II$ and
$(I,L)\in \K$ we have
\begin{align*}
 \dim A_{I}(\lambda,\gamma)
& \leq \dim \toba(M_{I}) \vert \DD_{m}\vert = 4^{|I|}2m\qquad  \text{ and }\\
\dim B_{I,L}(\lambda,\gamma,\theta,\mu)
& \leq \dim \toba(M_{I,L} )
\vert \DD_{m}\vert = 4^{|I|+|L|}2m.
\end{align*}\qed
\end{cor}

\subsection{Representation theory}
Let $H$ be a finite-dimensional pointed
Hopf algebra over $\k \DD_{m}$.
In this subsection we prove
using representation theory
that the
quadratic algebras defined
in Definitions \ref{liftings1} and \ref{def:liftings2}
are liftings of finite-dimensional
Nichols algebras over $ \k \DD_{m} $ for
all lifting data $(\lambda,\gamma)$
or $(\lambda,\gamma,\theta,\mu)$, and we
end the section with the proof of Theorem \ref{thm:B}.

\par
Let $H= A_{I}(\lambda,\gamma)$ with $I\in \II$ or
$H = B_{I,L }(\lambda,\gamma,\theta,\mu)$ with 
$(I,L) \in \K$. By definition,
the group $G(H)$ is a quotient of
$ \DD_{m} $;
in particular, any $H$-module is a $\k \DD_{m}$-module.
Denote by $\pi: \DD_{m} \twoheadrightarrow G(H)$ this quotient. The
following lemma is a key step to determine the dimension
of $H$.

\begin{lema}\label{lema:dimension}
Let $\rho: H \to \End (V)   $ be a representation
of $ H $ such that
\begin{itemize}
\item[$ (i) $] $\rho_{|G(H)}\circ \pi: \DD_{m} \to \End(V)$
is faithful and
\item[$ (ii) $]  if $H =A_{I}(\lambda,\gamma)$, then
$\rho(x_{p,q}) \notin\k\rho(G(H))$ for all $(p,q) \in I$ and
if $H =B_{I,L}(\lambda,\gamma,\theta,\mu)$, then
$\rho(x_{p,q}), \rho(z_{\ell})\notin\k\rho(G(H))$ for all
$(p,q) \in I$ and $\ell \in L$.
\end{itemize}
Then $\gr H=\toba(M)\#\k \DD_{m}$
and thus $\dim H=\dim\toba(M)|\DD_{m}|$.
\end{lema}

\pf
Let
$M =M_{I,L } $ and suppose
that
$H = B_{I,L }(\lambda,\gamma,\theta,\mu)$.
Since $G(H)$ is a quotient of $ \DD_{m} $, from $(i)$ it follows that
$G(H) \simeq \DD_{m}$. Thus $H$ is a pointed Hopf algebra
over $\DD_{m}$ and by Theorem \ref{teo:item-c},
$\gr H \simeq \toba(N) \# \k \DD_{m}$,
with $\toba(N)$ an exterior
algebra, see  Theorem \ref{teo:all-nichols-dm}.
Furthermore, by Lemma \ref{lema:quadratic-bound}
we have that $\dim \toba(N) \leq \dim \toba(M)$. But by $ (ii) $
the map $\varphi: M  \to H_{1}/H_{0}$, sending
$$
a_{p,q} \mapsto \bar{x}_{p,q},\quad
b_{p,q} \mapsto \bar{y}_{p,q},\quad
c_{\ell} \mapsto \bar{z}_{\ell},\quad
d_{\ell} \mapsto \bar{w}_{\ell},
$$
induces an injective map $\phi: M \to N$ in $ \yddm $
which implies that $\dim\toba(N)$ $\geq \dim \toba(M)$.
The proof for $H = A_{I}(\lambda,\gamma)$
is completely analogous.
\epf

\subsection{Proof of Theorem B}
Let $H$ be a finite-dimensional pointed Hopf algebra
with $ G(H)=\dm $. Then
$\gr H\simeq R\#\k \DD_{m}$
and by Theorem \ref{teo:item-c} the diagram $R$ is a Nichols algebra
$\toba(M)$ for some $M\in \yddm$ and consequently it is isomorphic to
one of the Hopf algebras of Theorem \ref{teo:all-nichols-dm}.

If $M\simeq
M_{I}$ with $I=\{(i,k)\}$ and $k\neq n$ or  
$M\simeq
M_{L}$ with $L\in \Lc$, 
then
$H \simeq  \toba(M)\#\k \DD_{m}$
by Corollary \ref{prop:sin-lifting-ik-1}. If  $M\simeq
M_{I}$ with $I\in \II$ and $|I|>0$, 
then by Corollary \ref{cor:rel-lifting-ikI} there exists
a lifting datum $(\lambda,\gamma)$ and an epimorphism of
Hopf algebras $A_{I}(\lambda,\gamma) \twoheadrightarrow H$.
Hence $\dim H \leq \dim A_{I}(\lambda,\gamma) \leq
\dim \toba (M_{I})|\DD_{m}|$. This implies that
$H\simeq A_{I}(\lambda,\gamma)$, since
$\dim H = \dim \gr H =  \dim \toba (M_{I})|\DD_{m}|$.
If $M\simeq
M_{I,L}$ with $(I,L)\in \Lc$, then using the same argument
as before with Corollary \ref{cor:rel-lifting-ikI-L} shows that $H\simeq
B_{I,L}(\lambda,\gamma,\theta,\mu)$.

For the converse, it is clear  that the
algebras listed in items $(a)$ and $(b)$
are liftings of Nichols algebras over $\dm$.
Thus, we need to show only
that the Hopf algebras $ A_{I}(\lambda,\gamma)$ and
$B_{I,L}(\lambda,\gamma, \theta,\mu)$
are liftings for all $I\in \II$, $(I,L)\in \K$ and for all
lifting data.

Assume first that $I\in \II$.
Following Lemma \ref{lema:dimension}, we give a representation for
$ A_{I}(\lambda,\gamma)$. Give $I$ an order and write $I =
((i_{1},k_{1}),\ldots ,(i_{r},k_{r}))$. Let $V$ be a vector space
with basis given by two families of vectors $\{u_{\alpha}\}$,
$\{v_{\alpha}\}$, indexed by all possible ordered monomials in the
letters $i_{s,1}, i_{s,2}$ for all $1\leq s\leq r$ such that each
letter appears at most once (set $u_{0},v_{0}$ if no letter
appears) and the order is given by $i_{s,p} < i_{t,p}$ for all
$p=1,2$ iff $s<t$, $i_{s,1}<i_{s,2}$ for all $1\leq s\leq r$ and
$i_{s,2}<i_{t,1}$ iff $s<t$; \textit{e.g.}
$v_{i_{1,1}i_{2,2}i_{3,1}i_{3,2}}$ is a basis element. In
particular, $\dim V = 2\dim \bigwedge M_{I}$.

For all $1\leq j< n$,
$V$ bears an $A_{I}(\lambda,\gamma)$-module
structure determined by
\begin{align*}
 g\cdot u_{0}  = v_{0},\quad h\cdot u_{0} = \omega^{j} u_{0},\quad
x_{i_{s},k_{s}}\cdot u_{0} = u_{i_{s,1}},
\quad
y_{i_{s},k_{s}}\cdot u_{0} = u_{i_{s,2}}\\
 (x_{i_{t},k_{t}}x_{i_{s},k_{s}})\cdot u_{0}  =
- u_{i_{s,1}i_{t,1}} + \delta_{k_{s},m-k_{t}}
\lambda_{i_{s},k_{s},i_{t},k_{t}}
(1-\omega^{j(i_{s}+i_{t})}) u_{0}
\text{ if } t>s,\\
(y_{i_{t},k_{t}}x_{i_{s},k_{s}})\cdot u_{0}  =
- u_{i_{s,1}i_{t,2}} + \delta_{s,t}\gamma_{i_{s},k_{s},i_{t},k_{t}}
(1-\omega^{j(i_{s}-i_{t})}) u_{0}
\text{ if } t\geq s.
\end{align*}
because of the defining relations of $A_{I}(\lambda,\gamma)$,
see Definition \ref{liftings1}. Hence,
$\rho_{|G(A_{I}(\lambda,\gamma))}\circ \pi: \DD_{m} \to \End(V)$
is faithful since $(\k\{u_{0},v_{0}\},\rho_{|G(A_{I}(\lambda,\gamma))})
 \simeq (\k^{2},\rho_{j})$, and
$\rho(x_{i_{s},k_{s}}) \notin \k\rho(G(A_{I}(\lambda,\gamma))$ by
definition. Then $\gr A_{I}(\lambda,\gamma) = \toba(M_{I})\# \k
\dm$ and $A_{I}(\lambda,\gamma)$ is a lifting. The proof for
$B_{I,L}(\lambda,\gamma, \theta ,\mu)$ is analogous.\qed

\subsubsection{Isomorphism classes}In this last subsection
we
study the isomorphism classes of the
families of
Hopf algebras given by Theorem \ref{teo:liftings-quadratic}.

Let $H$ be a finite-dimensional pointed
Hopf algebra over $ \dm $. Then
$H$ is isomorphic to a Hopf algebra listed in
Theorem \ref{teo:liftings-quadratic}; in particular, it is
a lifting of a finite-dimensional
Nichols algebra over $\dm$.

It is clear that two algebras from different families are not
isomorphic as Hopf algebras since their infinitesimal braidings
are not isomorphic as Yetter-Drinfeld modules.

Thus, we have to show that two different members in the same family
are not isomorphic.
By the argument above, if $I=\{(i,k)\}$, $I'=\{(p,q)\} \in \II$,
with $k,q\neq n$ and
$(i,k)\neq (p,q)$, then
$\toba(M_{I})\#\k \dm \not\simeq \toba(M_{I'})\#\k \dm$,
and if $L,L' \in \Lc$ with 
$L\neq L'$, then $\toba(M_{L})\#\k \dm \not\simeq \toba(M_{L'})\#\k \dm$.
We end the paper by showing the isomorphism classes of the
families of the items $(c)$ and $(d)$.

Observe that $\Z/m$ acts on $\II$ with the action on
each $I\in \II$ induced by
$$\ell\cdot (i_{s},k_{s})  = \begin{cases}\begin{matrix}
(\ell i_{s}, \ell^{-1}k_{s})\qquad
\text{ if }1\leq \ell i_{s} <n  \mod m,\\
(m-\ell i_{s}, \ell^{-1}k_{s})\qquad
\text{ if }n\leq \ell i_{s} \mod m.
\end{matrix}
\end{cases}$$

\begin{lema}\label{lema:iso-AI}
Let $I,I'\in\II$. Then  $A_{I}(\lambda,\gamma) \simeq
A_{I'}(\lambda',\gamma')$ if and only if
there exists $\ell \in \Z/m$, with $(\ell,m) = 1$ such that
 $\ell\cdot I = I'$,
and for all $(p,q),(i,k)\in I$,
\begin{eqnarray}
\begin{cases}
\lambda_{p,q,i,m-q}  = \lambda'_{\ell\cdot(p,q),\ell\cdot(i,m-q)},
\\
\gamma_{p,q,i,q}  = \gamma'_{\ell\cdot(p,q),\ell\cdot(i,q)},
\end{cases} \text{if }p\ell,\ i\ell  < n
\text{ or }n<p\ell,\ i\ell \mod m,\label{eq:iso-a-I}
\\
\text{ and }\qquad\begin{cases}
\delta_{q,m-k}\lambda_{p,q,i,k}
= \delta_{k,q}\gamma'_{\ell\cdot(p,q),\ell\cdot(i,k)},
\\
\delta_{q,k}\gamma_{p,q,i,q}
= \delta_{q,m-k}\lambda'_{\ell\cdot(p,q),\ell\cdot(i,q)},
\label{eq:iso-a-I-cruzado}
\end{cases}\text{ otherwise.}
\end{eqnarray}
In particular, $A_{I}(\lambda,\gamma) \simeq \toba(M_{I})\# \k \dm$ if
and only if $\lambda\equiv 0 \equiv \gamma$.
\end{lema}

\pf Suppose $\varphi: A_{I}(\lambda,\gamma) \to
A_{I'}(\lambda',\gamma')$ is a Hopf algebra isomorphism and denote
by $g,h,x_{i,k},y_{i,k}$ and $g',h',x'_{i,k},y'_{i,k}$ the
generators of $A_{I}(\lambda,\gamma)$ and
$A_{I'}(\lambda',\gamma')$, respectively. Since both must have the
same dimension, we have that $ |I| =|I'| $. Moreover, $\varphi(g)
= g'$, $\varphi(h) = h'^{\ell}$ for some $\ell \in \Z/m$ with
$(m,\ell) =1$, and $\varphi(x_{i,k})$, $\varphi(y_{i,k})$ are
$(h'^{i\ell},1) $-primitive and $(h'^{-i\ell},1) $-primitive in
$A_{I'}(\lambda',\gamma')$ for all $(i,k) \in I$, respectively.
Using that $\varphi(hx_{i,k}    h^{-1}) =
h'^{\ell}\varphi(x_{i,k})h'^{-\ell}$ we have that
$\varphi(x_{i,k}) = a_{i,k,i\ell,k\ell^{-1}}x'_{i\ell,k\ell^{-1}}$
if
 $i\ell<n$ and
$\varphi(x_{i,k}) = b_{i,k,m-i\ell,k\ell^{-1}}y'_{m-i\ell,k\ell^{-1}}$
if $n< i\ell<m$,
 for some $a_{i,k,p,q},b_{i,k,p,q}\in \k^{\times}$.
 In particular, this implies that $I' = \ell\cdot I$. Clearly,
we may assume that $a_{i,k,p,q},b_{i,k,p,q} = 1$. Denote
$\varphi_{\ell} = \varphi$.

Let $(p,q),(i,k) \in I$ and suppose that $\ell p,\ell i < n$. Then
applying $\varphi_{\ell}$ on both sides of \eqref{eq:x-x-a-I}
yields $$x'_{\ell p, \ell^{-1}q}x'_{\ell i, \ell^{-1}k}
+ x'_{\ell i, \ell^{-1}k}x'_{\ell p, \ell^{-1}q} =
\delta_{q,m-k} \lambda_{p,q,i,k}(1-h'^{\ell(p+i)}).
$$
But the left hand side equals $\delta_{\ell^{-1}q,m-\ell^{-1}k}
\lambda'_{\ell p,\ell^{-1}q,\ell i,\ell^{-1}k}(1-h'^{\ell p+\ell i}) $,
by the same relation in
$A_{I'}(\lambda',\gamma')$. Hence
$ \lambda_{p,q,i,m-q} = \lambda'_{\ell\cdot(p,q),\ell\cdot(i,m-q)}$.
On the other hand,
applying $\varphi_{\ell}$ on both sides of \eqref{eq:x-y-a-I}
yields $$x'_{\ell p, \ell^{-1}q}y'_{\ell i, \ell^{-1}k}
+ y'_{\ell i, \ell^{-1}k}x'_{\ell p, \ell^{-1}q} =
\delta_{q,k} \gamma_{p,q,i,k}(1-h'^{\ell(p-i)}).
$$
But the left hand side equals $\delta_{\ell^{-1}q,\ell^{-1}k}
\gamma'_{\ell p,\ell^{-1}q,\ell i,\ell^{-1}k}(1-h'^{\ell p-\ell i}) $,
by the same relation in
$A_{I'}(\lambda',\gamma')$. Hence
$ \gamma_{p,q,i,q} = \gamma'_{\ell\cdot(p,q),\ell\cdot(i,q)}$.
The proof for the remaining cases is completely analogous.

Assume now
there exists $\ell \in \Z/m$, with $(\ell,m) = 1$ such that
 $\ell\cdot I = I'$,
and equations \eqref{eq:iso-a-I}, \eqref{eq:iso-a-I-cruzado} hold
for all $(p,q),(i,k)\in I$. Then we may define an algebra morphism
by $\varphi_{\ell}(x_{p,q}) = x'_{\ell p, \ell^{-1}q}$,
$\varphi_{\ell}(y_{p,q}) = y'_{\ell p, \ell^{-1}q}$ if $\ell p <n$
and $\varphi_{\ell}(x_{p,q}) = y'_{m-\ell p, \ell^{-1}q}$,
$\varphi_{\ell}(y_{p,q}) = x'_{m-\ell p, \ell^{-1}q}$
if $n<\ell p \mod m$.
Equations \eqref{eq:iso-a-I}, \eqref{eq:iso-a-I-cruzado}
and
the fact that $\ell\cdot I = I'$,
ensure that $\varphi_{\ell}$ is a well-defined surjective
Hopf algebra map, which is indeed an isomorphism by Thm.
\ref{prop:quad-are-liftings}.
\epf

Note also that $\Z/m$ also acts on $\Lc$ with the action induced by
$$\ell\cdot r  = \begin{cases}
\begin{matrix}\ell^{-1}r\qquad
\text{ if }1\leq \ell^{-1} r <n \mod m,\\
m-\ell^{-1} r\qquad
\text{ if }n\leq \ell^{-1} r \mod m.
\end{matrix}
\end{cases}$$

The proof of the following lemma is completely
analogous to the proof of Lemma \ref{lema:iso-AI}.

\begin{lema}\label{lema:iso-BIL}
Let $(I,L)$, $(I',L')\in \K$. $B_{I,L}(\lambda,\gamma,\theta,\mu)$ $ \simeq
B_{I',L'}(\lambda',\gamma',\theta',\mu')$ if and only if
there exists $\ell \in \Z/m$ with $(\ell,m) = 1$ such that
 $\ell\cdot I = I'$, $\ell\cdot L = L'$, $\lambda, \lambda'$
 and $\gamma, \gamma'$ satisfy conditions \eqref{eq:iso-a-I} and
 \eqref{eq:iso-a-I-cruzado},
and for all $(p,q)\in I$, $r\in L$,
\begin{eqnarray}
\nonumber \begin{cases}
\delta_{q,m-r}\theta_{p,q,r}  =
\delta_{q,m-r} \theta'_{\ell\cdot(p,q),\ell\cdot r},
\\
\delta_{q,r}\mu_{p,q,r}  = \delta_{q,r}\mu'_{\ell\cdot(p,q),\ell\cdot r},
\end{cases} \text{if }p\ell,\ r\ell^{-1}  < n,\\
\nonumber \begin{cases}
\delta_{q,m-r}\theta_{p,q,r}  =
\delta_{q,r} \theta'_{\ell\cdot(p,q),\ell\cdot r},
\\
\delta_{q,r}\mu_{p,q,r}  = \delta_{q,m-r}\mu'_{\ell\cdot(p,q),\ell\cdot r},
\end{cases} \text{if } n < p\ell,\ r\ell^{-1}  \mod m,\\
\nonumber \begin{cases}
\delta_{q,m-r}\theta_{p,q,r}  =
\delta_{q,m-r} \mu'_{\ell\cdot(p,q),\ell\cdot r},
\\
\delta_{q,r}\mu_{p,q,r}  = \delta_{q,r}\theta'_{\ell\cdot(p,q),\ell\cdot r},
\end{cases} \text{if } p\ell <n< r\ell^{-1}  \mod m,\\
\nonumber \begin{cases}
\delta_{q,m-r}\theta_{p,q,r}  =
\delta_{q,r} \mu'_{\ell\cdot(p,q),\ell\cdot r},
\\
\delta_{q,r}\mu_{p,q,r}  = \delta_{q,m-r}\theta'_{\ell\cdot(p,q),\ell\cdot r},
\end{cases} \text{if } r\ell^{-1}< n < p\ell \mod m.
\end{eqnarray}\qed
\end{lema}

\subsection*{Acknowledgment}
The authors want to thank N. Andruskiewitsch for suggestions
and fruitful comments.

\end{document}